\documentclass[a4paper,11pt]{article}
\usepackage{amsmath,amsfonts,amsthm,dsfont,amssymb}
\usepackage[blocks]{authblk}
\usepackage{cite}
\usepackage{graphicx}
\usepackage{amssymb,array,arydshln}
\usepackage{amsfonts}
\usepackage{amsbsy}
\usepackage{cite}
\usepackage{amssymb}
\usepackage{amsmath}
\usepackage{graphics}
\newenvironment{keywords}{\noindent\textbf{Keywords:}}{}
\newenvironment{classification}{\noindent\textbf{AMS subject classifications.}}{}
\date{}
\newcommand{\email}[1]{\texttt{\small #1}}
\newtheorem{theorem}{Theorem}[section]
\newtheorem{remark}[theorem]{Remark}
\newtheorem{example}[theorem]{Example}
\newtheorem{lemma}[theorem]{Lemma}
\newtheorem{corollary}[theorem]{Corollary}
\newtheorem{definition}[theorem]{Definition}
\newtheorem{proposition}[theorem]{Proposition}

\newtheorem{conjecture}{Conjecture}
\newtheorem{question}{Question}
\newcommand{\bea}{\begin{eqnarray}}
\newcommand{\eea}{\end{eqnarray}}
\newcommand{\beq}{\begin{eqnarray*}}
\newcommand{\eeq}{\end{eqnarray*}}
\def \bd{\begin{definition}}
\def \ed{\end{definition}}
\def \bqu{\begin{question}}
\def \equ{\end{question}}
\def \bcc{\begin{conjecture}}
\def \ecc{\end{conjecture}}
\def \bt{\begin{theorem}}
\def \et{\end{theorem}}
\def \bl{\begin{lemma}}
\def \el{\end{lemma}}
\def \bc{\begin{corollary}}
\def \ec{\end{corollary}}
\def \be{\begin{equation}}
\def \ee{\end{equation}}
\def \ben{\begin{enumerate}}
\def \een{\end{enumerate}}
\def \ba{\begin{array}}
\def \ea{\end{array}}
\def \bp{\begin{proposition}}
\def \ep{\end{proposition}}
\def \bx{\begin{example}}
\def \ex{\end{example}}
\def \br{\begin{remark}}
\def \er{\end{remark}}
\def \bdsc{\begin{description}}
\def \edsc{\end{description}}

\def\1{1\!\!1}
\def\0{0\!\!0}

\begin{document}
	\title{On the distance \& distance (signless) Laplacian spectra of non-commuting graphs
	}

\author[1]{Jharna Kalita}
\author[2]{Somnath Paul\footnote{Corresponding  Author.}}
\affil[1]{Department of Applied Sciences\\ Tezpur University\\ Napaam-784028, Assam, India. \email{app21104@tezu.ac.in}}
\affil[2]{Department of Applied Sciences\\ Tezpur University\\ Napaam-784028, Assam, India. \email{som@tezu.ernet.in}}

\pagestyle{myheadings} \markboth{J. Kalita \& S. Paul}{On the distance \& distance (signless) Laplacian spectra of non-commuting graphs}
	\maketitle
	
\begin{abstract}
Let $Z(G)$ be the centre of a finite non-abelian group $G.$ The non-commuting graph of $G$  is a simple undirected graph with vertex set $G\setminus Z(G),$ and two vertices $u$ and $v$ are adjacent if and only if $uv\ne vu.$ In this paper, we investigate the distance, distance (signless) Laplacian spectra of non-commuting graphs of some classes of finite non-abelian groups, and obtain some conditions on a group so that the non-commuting graph is distance, distance (signless) Laplacian integral.
\end{abstract}


\section{Introduction}

Let $G$ be a finite non-ableian group and $Z(G)$ be its centre. The non-commuting graph of $G,$ denoted by $\Gamma_G,$ is a simple undirected graph with vertex set $G\setminus Z(G),$ and two vertices $u$ and $v$ are adjacent if and only if $uv\ne vu.$ In literature, various studies have been done discussing several aspects of non-commuting graphs of different finite groups (see \cite{1,4,8,23}).

For a simple graph $H,$ the \textit{distance} between two of its vertices $u$ and $v$ is defined as the length of a shortest path between $u$ and $v$ in it, and is denoted by $d_{uv}.$ The \textit{distance matrix} of $H$ is denoted by $D(H),$ with the $(u,v)$-entry equal to $d_{uv}.$

The \textit{transmission} of a vertex $v$ is defined to be the sum of the distances from $v$ to all other vertices in the graph. Let $Tr(H)$ denote the diagonal matrix with $i$-th diagonal entry being the  transmission of the $i$-th vertex in $H.$ In \cite{AH2}, M. Aouchiche and P. Hansen have defined the Laplacian and the signless Laplacian for the distance matrix of a connected graph, which are analogous to the Laplacian matrix and signless Laplacian matrix for adjacency matrix. The matrix $ D^L(H)=Tr(H)-D(H) $ is called the \textit{distance Laplacian matrix} of $H$, while the matrix  $D^{Q}(H)=Tr(H)+D(H)$ is called the \textit{distance signless Laplacian matrix of $H$}. For a beautiful survey of distance matrix, the reader can see \cite{AH1}.

If $M$ is a symmetric matrix, then the characteristic polynomial of $M$ has only real zeroes. We will represent this family of eigenvalues (known as the \textit{spectrum}) as $$\sigma_M=\left(
                  \begin{array}{cccc}
                    \mu_1 & \mu_2 &\cdots&\mu_p\\
                    m_1 & m_2 &\cdots & m_p
                  \end{array}
                \right),
$$ where $\mu_1,\mu_2,\ldots,\mu_p$ are the distinct eigenvalues of $M$ and $m_1, m_2,\ldots,m_p$ are the corresponding multiplicities. Since each of $D(H),~D^L(H)$ and $D^Q(H)$ is symmetric, we will refer the corresponding spectrum as the \textit{distance, distance Laplacian} and \textit{distance signless Laplacian spectrum}, respectively.

A group $G$ is called a \textit{CA-group} if the centralizer $C_G(x)$ for every non central element $x$ is abelian. For any $m\ge3,$ the authors of \cite{rkn}  used the concept of CA group and its centralizers for computing the Laplacian spectrum of the non-commuting graph of $D_{2m}.$ Not only this, they have also used the method to compute the Laplacian spectrum of the non-commuting graphs of several well-known families of finite non-abelian groups such as the quasidihedral groups, generalized quaternion groups, some projective special linear groups, general linear groups etc. Though the distance spectra of non-commuting graph of $D_{2m}$ is discussed in \cite{ban}, but not many results are available in the literature discussing the distance spectra of non-commuting graphs of other groups. Here, we tried to fill the gap and obtain  the distance, distance (signless) Laplacian spectra of non-commuting graphs of some classes of CA  groups, namely  the generalized quaternion group, the quasidihedral Group, $U_{6n}$ and the metacyclic group.

A graph is called distance (respectively distance (signless) Laplacian) integral if the distance (respectively distance (signless) Laplacian) spectrum consists entirely of integers.  In this article, we consider four CA groups and obtain some conditions on them, so that their non-commuting graph is distance, distance (signless) Laplacian integral.

\section{Preliminaries}

Let $G$ be a finite non-abelian CA-group and $X_1,X_2,\ldots,X_n$ are the distinct centralizers of non-central elements of $G,$ then it is observed in \cite{rkn} that $\Gamma_G$ is the complete multipartite graph $K_{m_1,m_2,\ldots,m_n},$ where $m_i=|X_i|-|Z(G)|,$ for $i=1,2,\ldots,n.$

Consider $Q_{4n} =<x,y : x^{2n} = 1, x^n = y^2, yx=x^{-1}y>,$ be the generalized quaternion group of order $4n,$  where $n\ge2.$ Then $Q_{4n}$ can be written as $A\cup B,$ where $A=\{1,x,x^2,\ldots,x^{2n-1}\}$ and $B=\{y,xy,x^2y,\ldots,x^{2n-1}y\},$ where each element of $B$ is of order 4. It can be easily verified that $Z(Q_{4n})=\{1,x^n\},$ and thus $|Z(Q_{4n})|=2.$ Therefore, from the definition of non-commuting graph, it follows that $\Gamma_{Q_{4n}}$ is a graph with vertex set $(A\setminus Z(Q_{4n}))\cup B,$ where no two vertices of $A\setminus Z(Q_{4n})$ are adjacent to each other.  Also every vertex of $A\setminus Z(Q_{4n})$ is adjacent to every vertex of $B.$ Moreover, it can also be seen that two vertices $x^iy$ and $x^jy$ of $B$ are adjacent \textit{iff} $|i-j|=n.$ Therefore, $\Gamma_{Q_{4n}}$ is the complete $(n+1)$-partite graph $K_{2n-2,\underbrace{2,2,\ldots,2}_n}.$

We also note that for any $z\in Z(Q_{4n}),$ and $1\le i\le 2n-1,$
\begin{eqnarray*}
  C_{Q_{4n}}(x)
  &=&C_{Q_{4n}}(x^iz)\\ &=& Z(Q_{4n})\cup x Z(Q_{4n})\cup x^2Z(Q_{4n})\cup\ldots\cup x^{n-1}Z(Q_{4n})\\
    &=& \{1,x^n\} \cup x \{1,x^n\}\cup x^2\{1,x^n\}\cup\ldots\cup x^{n-1}\{1,x^n\}\\
    &=& \{1,x^n\} \cup \{x,x^{n+1}\}\cup \{x^2,x^{n+2}\}\cup\ldots\cup \{x^{n-1},x^{2n-1}\} \\
    &=& \{1,x,x^2,\ldots,x^{2n-1}\}.
\end{eqnarray*}
Moreover, for $1\le j\le n,$
\begin{eqnarray*}
  C_{Q_{4n}}(x^jy)=C_{Q_{4n}}(yx^jz) &=& Z(Q_{4n})\cup yx^j Z(Q_{4n}) \\
    &=& \{1,x^n\} \cup yx^j \{1,x^n\}\\
    &=& \{1,x^n\} \cup \{yx^j,yx^{n+j}\}.
\end{eqnarray*}

Thus the structure of $\Gamma_{Q_{4n}}$ given in \cite{rkn}, using the concept of centralizers of $Q_{4n}$ and showing that $Q_{4n}$ is a CA group matches with the structure discussed in the preceding paragraph.

The following result gives the distance characteristics polynomial for $K_{n_1,n_2,\ldots,n_k}$ and will be useful to derive some of our main results.
\bl{\rm \cite{disc}}\label{lem1} The distance characteristics polynomial of the complete multipartite graph $K_{n_1,n_2,\ldots,n_k}$ is
\begin{equation}\label{e1}
 P_D(\lambda)=(\lambda+2)^{n-k}\left[\prod_{i=1}^k(\lambda-n_i+2)-\sum_{i=1}^k n_i\prod_{j=1,j\neq i}^k(\lambda-n_j+2)\right].
\end{equation} \el
\section{Spectra of $\Gamma_{Q_{4n}}$}\label{sec1}

In this section, we consider the non commuting graph of $Q_{4n}$ and obtain the distance, distance Laplacian and distance signless Laplacian spectra of it. As it is already observed in the previous section,  $\Gamma_{Q_{4n}}=K_{2n-2,\underbrace{2,2,\ldots,2}_n}$.
\subsection{Distance spectrum of $\Gamma_{Q_{4n}}$}
Since $\Gamma_{Q_{4n}}$ is a complete $(n+1)$-partite graph of order $4n-2$ with $n$ partitions of cardinality 2 and one partition of cardinality $2n-2,$ by equation (\ref{e1}) we get,
\begin{eqnarray*}
  &&P_{D(\Gamma_{Q_{4n}})}(\lambda) \\
    &=& (\lambda+2)^{(3n-3)}\left[(\lambda-2n+4)\lambda^n-(2n-2)\lambda^n-2n\lambda^{(n-1)}(\lambda-2n+4)\right]  \\
    &=&(\lambda+2)^{(3n-3)}\left[\lambda^{n+1}-(6n-6)\lambda^n+4n(n-2)\lambda^{n-1}\right]\\
    &=&(\lambda+2)^{(3n-3)}\lambda^{n-1}[\lambda^2-6(n-1)\lambda+4n(n-2)]
\end{eqnarray*}

Therefore, eigenvalues of $D(\Gamma_{Q_{4n}})$ are $-2$ with multiplicity $(3n-3),$~$0$ with multiplicity $(n-1),$ and two roots of the equation $\lambda^2-6(n-1)\lambda+4n(n-2)=0.$ Hence, we have the following theorem.

\bt\label{dt1} Let $\Gamma_{Q_{4n}}$ be the non commuting graph of the generalized quaternion group of order $4n.$ Then
\begin{enumerate}
\item[{\rm(a)}] $-2\in\sigma(D(\Gamma_{Q_{4n}}))$ with multiplicity $3n-3;$
\item[{\rm(b)}] $0\in\sigma(D(\Gamma_{Q_{4n}}))$ with multiplicity $n-1;$
\item[{\rm(c)}] $3(n-1)\pm\sqrt{5n^2-10n+9}\in\sigma(D(\Gamma_{Q_{4n}})),$ each with multiplicity 1.
\end{enumerate}\et
Therefore, by Theorem~\ref{dt1}, it follows that $\Gamma_{Q_{4n}}$ is distance integral if  $5n^2-10n+9$ is a perfect square.
 \subsection{Distance Laplacian spectrum of $\Gamma_{Q_{4n}}$}
Let $\1_n$ (resp. $\0_n$) denote the $n\times 1$ vector with each entry 1 (resp. 0). Also, let $J_n$ denote the matrix of order $n$ with all entries equal to 1 (we will write $J$ if the order is clear from the context). It is well known that the distance Laplacian matrix for any graph is positive semidefinite with 0 being the smallest eigenvalue with multiplicity 1 and $\1_n$ is the corresponding eigenvector. The following theorem describes the distance Laplacian spectrum of $\Gamma_{Q_{4n}}.$

\bt\label{dlt1} Let $\Gamma_{Q_{4n}}$ be the non commuting graph of the generalized quaternion group of order $4n.$ Then
\begin{enumerate}
\item[{\rm(a)}] $0\in\sigma(D^L(\Gamma_{Q_{4n}}))$ with multiplicity $1;$
\item[{\rm(b)}] $4n-2,4n\in\sigma(D^L(\Gamma_{Q_{4n}})),$ each with multiplicity $n;$
\item[{\rm(c)}] $6n-4\in\sigma(D^L(\Gamma_{Q_{4n}}))$  with multiplicity $2n-3.$
\end{enumerate}\et
 {\bf Proof.} With a suitable labeling of the vertices, the distance Laplacian matrix for $\Gamma_{Q_{4n}}$ can be written as
$$D^L(\Gamma_{Q_{4n}})=\left[
                        \begin{array}{c|c|c|c}
                          (6n-4)I-2J & -J & \ldots & -J \\
                          \hline
                          -J & 4nI-2J & \ldots & -J \\
                          \hline
                          \ldots & \ldots & \ldots & \ldots \\
                          \hline
                          -J & -J & \ldots & 4nI-2J \\
                        \end{array}
                      \right]$$

Obviously, 0 is an eigenvalue of $D^L(\Gamma_{Q_{4n}})$ with multiplicity 1 and $\1_{4n-2}$ is the corresponding eigenvector.

Also, $D^L(\Gamma_{Q_{4n}})\left(
                             \begin{array}{c}
                               -\1_{2n-2} \\
                               \hline
                               (n-1)\1_2 \\
                               \hline
                               \0_{2n-2} \\
                             \end{array}
                           \right)
                           =(4n-2)\left(
                             \begin{array}{c}
                               -\1_{2n-2} \\
                               \hline
                               (n-1)\1_2 \\
                               \hline
                               \0_{2n-2} \\
                             \end{array}
                           \right).$

Therefore, $4n-2$ is an eigenvalue of $D^L(\Gamma_{Q_{4n}}),$ and in this way we can construct the following set $S_1$ of $n$ independent eigenvectors corresponding to $4n-2$.
$$S_1=\left\{\left(
         \begin{array}{c}
           -\1_{2n-2} \\
           \hline
           (n-1)\1_2 \\
           \hline
           \0_{2n-2} \\
         \end{array}
       \right),\left(
                 \begin{array}{c}
                   -\1_{2n-2} \\
                   \hline
                   \0_2 \\
                   \hline
                    (n-1)\1_2\\
                    \hline
                  \0_{2n-4} \\
                 \end{array}
               \right),\ldots,\left(
                                \begin{array}{c}
                                  -\1_{2n-2} \\
                                  \hline
                                   \0_{2n-2}\\
                                   \hline
                                  (n-1)\1_2 \\
                                \end{array}
                              \right)\right\}.$$

 Similarly, $D^L(\Gamma_{Q_{4n}})\left(
                                   \begin{array}{c}
                                   \0_{2n-2} \\
                                   \hline
                                     -1 \\
                                     1 \\
                                     \hline
                                     \0_{2n-2} \\
                                   \end{array}
                                 \right)=4n\left(
                                   \begin{array}{c}
                                   \0_{2n-2} \\
                                   \hline
                                     -1 \\
                                     1 \\
                                     \hline
                                     \0_{2n-2} \\
                                   \end{array}
                                 \right),$
shows that $4n$ is an eigenvalue of $D^L(\Gamma_{Q_{4n}})$. In this way we can construct the following set $S_2$ of $n$ independent eigenvectors corresponding to $4n$.
                                           $$S_2=\left\{\left(
                                                   \begin{array}{c}
                                                      \0_{2n-2}\\
                                                      \hline
                                                     -1 \\
                                                     1 \\
                                                     \hline
                                                     \0_{2n-2} \\
                                                   \end{array}
                                                 \right),\left(
                                                           \begin{array}{c}
                                                             \0_{2n-2} \\
                                                             \hline
                                                             \0_2 \\
                                                             \hline
                                                             -1 \\
                                                             1 \\
                                                             \hline
                                                             \0_{2n-4} \\
                                                           \end{array}
                                                         \right),\ldots,\left(
                                                                          \begin{array}{c}
                                                                            \0_{2n-2} \\
                                                                            \hline
                                                                             \0_{2n-2}\\
                                                                             \hline
                                                                             -1\\
                                                                            1\\
                                                                          \end{array}
                                                                        \right)\right\}.$$

 Finally, $D^L(\Gamma_{Q_{4n}})\left(
                                  \begin{array}{c}
                                    -1 \\
                                     1 \\
                                    \0_{2n-4} \\
                                    \hline
                                    \0_{2n} \\
                                  \end{array}
                                \right)=(6n-4)\left(
                                  \begin{array}{c}
                                    -1 \\
                                     1 \\
                                    \0_{2n-4} \\
                                    \hline
                                    \0_{2n} \\
                                  \end{array}
                                \right).$

Thus, $6n-4$ is an eigenvalue of $D^L(\Gamma_{Q_{4n}}),$ and in this way we can construct the following set $S_3$ of $2n-3$ independent eigenvectors corresponding to $6n-4$.
$$S_3=\left\{\left(
        \begin{array}{c}
          -1 \\
          1 \\
          \0_{2n-4} \\
          \hline
          \0_{2n} \\
        \end{array}
      \right),\left(
                \begin{array}{c}
                  -1 \\
                  0 \\
                  1 \\
                  \0_{2n-5} \\
                  \hline
                  \0_{2n} \\
                \end{array}
              \right),\ldots,\left(
                               \begin{array}{c}
                                 -1 \\
                                 \0_{2n-4} \\
                                 1 \\
                                 \hline
                                 \0_{2n} \\
                               \end{array}
                             \right)\right\}.$$

It can be seen that $\1_{4n-2}\cup S_1\cup S_2\cup S_3$ is a set of mutually orthogonal vectors. Since the order of $\Gamma_{Q_{4n}}$ is $4n-2,$ the result follows.\qed

Hence, from Theorem~\ref{dlt1} it follows that $\Gamma_{Q_{4n}}$ is distance Laplacian integral for all $n.$

\subsection{Distance signless Laplacian spectrum of $\Gamma_{Q_{4n}}$}

The following theorem describes the distance signless Laplacian spectrum of $\Gamma_{Q_{4n}}.$

\bt\label{dqt1} Let $\Gamma_{Q_{4n}}$ be the non commuting graph of the generalized quaternion group of order $4n.$ Then
\begin{enumerate}
\item[{\rm(a)}] $4n-4\in\sigma(D^Q(\Gamma_{Q_{4n}}))$ with multiplicity $n;$
\item[{\rm(b)}] $6n-8 \in\sigma(D^Q(\Gamma_{Q_{4n}}))$ with multiplicity $2n-3;$
\item[{\rm(c)}] $4n-2\in\sigma(D^Q(\Gamma_{Q_{4n}}))$  with multiplicity $n-1.$
\item[{\rm (d)}] $t_k(2n-2)+(6n-2)\in\sigma(D^Q(\Gamma_{Q_{4n}}))$ with multiplicity 1, where $t_k$ is a root of the equation $x^2(2n-2)+(-4n+10)x-2n=0,$ for each $k=1,2.$
\end{enumerate}\et
{\bf Proof.} With a suitable labeling of the vertices, the distance signless Laplacian matrix for $\Gamma_{Q_{4n}}$ can be written as
$$D^Q(\Gamma_{Q_{4n}})=\left[
                        \begin{array}{c|c|c|c}
                          (6n-8)I+2J & J & \ldots & J \\
                          \hline
                          J & (4n-4)I+2J & \ldots & J \\
                          \hline
                          \ldots & \ldots & \ldots & \ldots \\
                          \hline
                          J & J & \ldots & (4n-4)I+2J \\
                        \end{array}
                      \right].$$

Now, $D^Q(\Gamma_{Q_{4n}})\left(
           \begin{array}{c}
             \0_{2n-2} \\
             \hline
             -1 \\
             1\\
             \hline
             \0_{2n-2} \\
           \end{array}
         \right)=(4n-4)\left(
           \begin{array}{c}
             \0_{2n-2} \\
             \hline
             -1 \\
             1\\
             \hline
             \0_{2n-2} \\
           \end{array}
         \right).$

Therefore, $(4n-4)$ is an eigenvalue of $D^Q(\Gamma_{Q_{4n}}),$ and in this way we can construct the following set $S_1$ of $n$ independent eigenvectors corresponding to $(4n-4)$.
$$S_1= \left\{\left(
          \begin{array}{c}
            \0_{2n-2} \\
            \hline
            -1 \\
            1 \\
            \hline
            \0_{2n-2} \\
          \end{array}
        \right),\left(
                  \begin{array}{c}
                    \0_{2n-2} \\
                    \hline
                    \0_2 \\
                    \hline
                    -1 \\
                    1 \\
                    \hline
                    \0_{2n-4} \\
                  \end{array}
                \right),\ldots,\left(
                                 \begin{array}{c}
                                  \0_{2n-2} \\
                                   \hline
                                   \0_{2n-2} \\
                                   \hline
                                   -1 \\
                                   1 \\
                                 \end{array}
                               \right)\right\}.$$

Similarly, $D^Q(\Gamma_{Q_{4n}})\left(
                 \begin{array}{c}
                   -1 \\
                   1 \\
                   \0_{2n-4} \\
                   \hline
                   \0_{2n} \\
                 \end{array}
               \right)=(6n-8)\left(
                               \begin{array}{c}
                                 -1 \\
                                 1 \\
                                 \0_{2n-4} \\
                                 \hline
                                 \0_{2n} \\
                               \end{array}
                             \right),$
shows that  ${(6n-8)}$ is an eigenvalue of $D^Q(\Gamma_{Q_{4n}}).$ In this way we can construct the following set $S_2$ of ${(2n-3)}$ independent eigenvectors corresponding to ${(6n-8)}.$

$$ S_2=\left\{\left(
         \begin{array}{c}
           -1 \\
           1 \\
           \0_{2n-4} \\
           \hline
           \0_{2n} \\
         \end{array}
       \right),\left(
                 \begin{array}{c}
                   -1 \\
                   0 \\
                   1 \\
                   \0_{2n-5} \\
                   \hline
                   \0_{2n} \\
                 \end{array}
               \right), \ldots, \left(
                                  \begin{array}{c}
                                    -1 \\
                                    \0_{2n-4} \\
                                    1 \\
                                    \hline
                                    \0_{2n} \\
                                  \end{array}
                                \right)\right\}.$$

Moreover, $D^Q(\Gamma_{Q_{4n}})\left(
                 \begin{array}{c}
                   \0_{2n-2} \\
                   \hline
                   -\1_2 \\
                                      \hline
                   \1_2 \\

                   \hline
                   \0_{2n-4} \\
                 \end{array}
               \right)=(4n-2)\left(
                 \begin{array}{c}
                   \0_{2n-2} \\
                   \hline
                   -\1_2 \\
                                      \hline
                   \1_2 \\

                   \hline
                   \0_{2n-4} \\
                 \end{array}
               \right).$
So, $(4n-2)$ is an eigenvalue of $D^Q(\Gamma_{Q_{4n}}),$ and in this way we can construct the following set $S_3$ of $(n-1)$ independent eigenvectors corresponding to $(4n-2)$.
$$S_3= \left\{\left(
        \begin{array}{c}
          \0_{2n-2} \\
          \hline
          -\1_2\\
          \hline
          \1_2 \\
                    \hline
          \0_{2n-4} \\
        \end{array}
      \right),\left(
                \begin{array}{c}
                  \0_{2n-2} \\
                  \hline
                  -\1_2 \\
                                    \hline
                  \0_2 \\
                                    \hline
                  \1_2 \\

                  \hline
                  \0_{2n-6} \\
                \end{array}
              \right), \ldots,\left(
                                \begin{array}{c}
                                  \0_{2n-2} \\
                                  \hline
                                  -\1_2 \\
                                                                    \hline
                                  \0_{2n-4} \\
                                  \hline
                                  \1_2 \\
                                                                  \end{array}
                              \right)\right\}.$$

Suppose $t$ be a scalar such that $\mu$ is an eigenvalue of $D^Q(\Gamma_{Q_{4n}})$ with eigenvector $ \left[
      \begin{array}{c}
        t\1_{2n-2}\\
               \hline

        \1_{2n} \\
          \end{array}
    \right].$
In that case, we have
\begin{eqnarray}
 \label{eq1} (10n-12)t+2n &=& \mu t \\
 \label{eq2} (2n-2)t+6n-2 &=& \mu
\end{eqnarray}

Eliminating $\mu$  from (\ref{eq1}) and (\ref{eq2}), we get $t$ is a root of the equation $$x^2(2n-2)+(-4n+10)x-2n=0.$$

Therefore, $(2n-2)t+(6n-2)$ is an eigenvalue of $D(\Gamma_{Q_{4n}}).$

It can be seen that $S_1\cup S_2\cup S_3\cup\left\{\left[
      \begin{array}{c}
        t_1\1_{2n-2}\\
               \hline

        \1_{2n} \\
          \end{array}
    \right], \left[
      \begin{array}{c}
        t_2\1_{2n-2}\\
               \hline

        \1_{2n} \\
          \end{array}
    \right]\right\}$   is a set of mutually orthogonal vectors, where $t_k$ is a root of $x^2(2n-2)+(-4n+10)x-2n=0,$ and $k=1,2.$ Since the order of $\Gamma_{Q_{4n}}$ is $4n-2,$ the result follows.\qed

Hence, from Theorem~\ref{dqt1}, it follows that the distance signless Laplacian spectrum of $\Gamma_{Q_{4n}}$ is integral for those values of $n$ which gives integral values to $t_k,$ for $k=1,2.$

\section{Spectra of $\Gamma_{QD_{2^n}}$}

In this section, we consider the non commuting graph of the Quasidihedral group $QD_{2^n}=<a,b:a^{2^{n-1}}=b^2=1,bab^{-1}=a^{2^{n-2}-1}>,$ where $n\ge4,$ and obtain the distance, distance Laplacian and distance signless Laplacian spectra of it. We will mention the statements of the results only, proofs of which are avoided as these are analogous to the proofs of the results obtained in section~\ref{sec1}.

It is observed in \cite{rkn} that  $\Gamma_{QD_{2^n}}$ is the complete $(2^{n-2}+1)$-partite graph $K_{2^{(n-1)}-2,\underbrace{2,2,\ldots,2}_{2^{n-2}}}.$
\bt\label{dg2t1} Let $\Gamma_{QD_{2^n}}$ be the non commuting graph of the Quasidihedral group of order $2^n.$ Then
\begin{enumerate}
\item[{\rm(a)}] $-2\in\sigma(D(\Gamma_{QD_{2^n}}))$ with multiplicity $3\times2^{n-2}-3;$
\item[{\rm(b)}] $0\in\sigma(D(\Gamma_{QD_{2^n}}))$ with multiplicity $2^{n-2}-1;$
\item[{\rm(c)}] $3(2^{n-2}-1)\pm\sqrt{5(2^{n-2})^2-10(2^{n-2})+9},$ each with multiplicity 1.
\end{enumerate}\et

Therefore, by Theorem~\ref{dg2t1}, it follows that $\Gamma_{QD_{2^n}}$ is distance integral if  $5(2^{n-2})^2-10(2^{n-2})+9$ is a perfect square.

\bt\label{dg2lt1} Let $\Gamma_{QD_{2^n}}$ be the non commuting graph of the Quasidihedral group of order $2^n.$ Then
\begin{enumerate}
\item[{\rm(a)}] $0\in\sigma(D^L(\Gamma_{QD_{2^n}}))$ with multiplicity $1;$
\item[{\rm(b)}] $2^n-2,2^n\in\sigma(D^L(\Gamma_{QD_{2^n}})),$ each with multiplicity $2^{n-2};$
\item[{\rm(c)}] $2^n+2^{n-1}-4\in\sigma(D^L(\Gamma_{QD_{2^n}}))$  with multiplicity $2^{n-1}-3.$
\end{enumerate}\et

Hence, from Theorem~\ref{dg2lt1} it follows that $\Gamma_{QD_{2^n}}$ is distance Laplacian integral for all $n.$
\bt\label{dg2qt1} Let $\Gamma_{QD_{2^n}}$ be the non commuting graph of the Quasidihedral group of order $2^n.$ Then
\begin{enumerate}
\item[{\rm(a)}] $2^n-4\in\sigma(D^Q(\Gamma_{QD_{2^n}}))$ with multiplicity $2^{n-2};$
\item[{\rm(b)}] $2^n-2 \in\sigma(D^Q(\Gamma_{QD_{2^n}}))$ with multiplicity $2^{n-2}-1;$
\item[{\rm(c)}] $2^n+2^{n-1}-8\in\sigma(D^Q(\Gamma_{QD_{2^n}}))$  with multiplicity $2^{n-1}-3.$
\item[{\rm (d)}] $t_k(2^{n-1}-2)+3(2^{n-1}-2) \in\sigma(D^Q(\Gamma_{QD_{4n}}))$ with multiplicity 1, where $t_k$ is a root of the equation $x^2(2^{n-1}-2)-(2^n-10)x-2^{n-1}=0,$ for each $k=1,2.$
\end{enumerate}\et

Hence, from Theorem~\ref{dg2qt1}, it follows that the distance signless Laplacian spectrum of $\Gamma_{QD_{2^n}}$ is integral for those values of $n$ which gives integral values to $t_k,$ for $k=1,2.$

\section{Spectra of $\Gamma_{U_{6n}}$}

In this section, we consider the non commuting graph of the  group $U_{6n}=<a,b: a^{2n}=b^3=1,a^{-1}ba=b^{-1}>,$ and obtain the distance, distance Laplacian and distance signless Laplacian spectra of it. We will mention the statements of the results only, proofs of which are avoided as these are analogous to the proofs of the results obtained in section~\ref{sec1}.

It can be easily observed that  $\Gamma_{U_{6n}}$ is the complete 4-partite graph $K_{2n,n,n,n}.$
\bt\label{dg3t1} Let $\Gamma_{U_{6n}}$ be the non commuting graph of the group $U_{6n}$ of order $6n.$ Then
\begin{enumerate}
\item[{\rm(a)}] $-2\in\sigma(D(\Gamma_{U_{6n}}))$ with multiplicity $5n-4;$
\item[{\rm(b)}] $n-2\in\sigma(D(\Gamma_{U_{6n}}))$ with multiplicity $2;$
\item[{\rm(c)}] $(4n-2)\pm n\sqrt{6}\in\sigma(D(\Gamma_{U_{6n}})),$ each with multiplicity 1.

\end{enumerate}\et

Therefore, by Theorem~\ref{dg3t1}, it follows that $\Gamma_{U_{6n}}$ is never distance integral.

\bt\label{dg3lt1} Let $\Gamma_{U_{6n}}$ be the non commuting graph of the group $U_{6n}$ of order $6n.$ Then
\begin{enumerate}
\item[{\rm(a)}] $0\in\sigma(D^L(\Gamma_{U_{6n}}))$ with multiplicity $1;$
\item[{\rm(b)}] $5n\in\sigma(D^L(\Gamma_{U_{6n}}))$  with multiplicity $3;$
\item[{\rm(c)}] $6n\in\sigma(D^L(\Gamma_{U_{6n}}))$  with multiplicity $3(n-1);$
\item[{\rm(d)}] $7n\in\sigma(D^L(\Gamma_{U_{6n}}))$  with multiplicity $2n-1.$
\end{enumerate}\et

Hence, from Theorem~\ref{dg3lt1} it follows that $\Gamma_{U_{6n}}$ is distance Laplacian integral for all $n.$
\bt\label{dg3qt1} Let $\Gamma_{U_{6n}}$ be the non commuting graph of the  group $U_{6n}$ of order $6n.$ Then
\begin{enumerate}
\item[{\rm(a)}] $6n-4\in\sigma(D^Q(\Gamma_{U_{6n}}))$ with multiplicity $3(n-1);$
\item[{\rm(b)}] $7n-4 \in\sigma(D^Q(\Gamma_{U_{6n}}))$ with multiplicity $2n+1;$
\item[{\rm (c)}] $8n-4 \in\sigma(D^Q(\Gamma_{U_{6n}}))$ with multiplicity $1;$
\item[{\rm (d)}] $13n-4 \in\sigma(D^Q(\Gamma_{U_{6n}}))$ with multiplicity $1.$
\end{enumerate}\et

Hence, from Theorem~\ref{dg3qt1}, it follows that the distance signless Laplacian spectrum of $\Gamma_{U_{6n}}$ is always integral.

\section{Spectra of $\Gamma_{M_{2mn}}$}

In this section, we consider the non commuting graph of the  Metacyclic group $M_{2mn}=<a,b: a^m=b^{2n}=1,bab^{-1}=a^{-1}>,$ and obtain the distance, distance Laplacian and distance signless Laplacian spectra of it, where $m>2.$  We will mention the statements of the results only, proofs of which are avoided as these are analogous to the proofs of the results obtained in section~\ref{sec1}.

It is observed in \cite{rkn} that  $\Gamma_{M_{2mn}}$ is the complete $(m+1)$-partite graph $K_{(m-1)n,\underbrace{n,n,\ldots,n}_m},$ when $m$ is odd, and complete $(\frac{m}{2}+1)$-partite graph $K_{(\frac{m}{2}-1)2n,\underbrace{2n,2n\ldots,2n}_{m/2}},$ when $m$ is even.

\bt\label{dg4ot1} Let $\Gamma_{M_{2mn}}$ be the non commuting graph of the Metacyclic group of order $2mn,$ where $m>2$ is odd. Then
\begin{enumerate}
\item[{\rm(a)}] $-2\in\sigma(D(\Gamma_{M_{2mn}}))$ with multiplicity $2mn-(m+n)-1;$
\item[{\rm(b)}] $n-2\in\sigma(D(\Gamma_{M_{2mn}}))$ with multiplicity $m-1;$
\item[{\rm(c)}] $\frac{-(4+n-3mn) \pm n\sqrt{5m^2-10m+9}}{2}\in\sigma(D(\Gamma_{M_{2mn}})),$ each with multiplicity 1.

\end{enumerate}\et
\bt\label{dg4et1} Let $\Gamma_{M_{2mn}}$ be the non commuting graph of the Metacyclic group of order $2mn,$ where $m>2$ is even. Then
\begin{enumerate}
\item[{\rm(a)}] $-2\in\sigma(D(\Gamma_{M_{2mn}}))$ with multiplicity $2n(m-1)-\frac{m}{2}-1;$
\item[{\rm(b)}] $2n-2\in\sigma(D(\Gamma_{M_{2mn}}))$ with multiplicity $\frac{m}{2}-1;$
\item[{\rm(c)}] $\frac{-(2n-3mn+4) \pm n\sqrt{5m^2-20m+36}}{2}\in\sigma(D(\Gamma_{M_{2mn}})),$ each with multiplicity 1.

\end{enumerate}\et

\bt\label{dg4olt1} Let $\Gamma_{M_{2mn}}$ be the non commuting graph of the Metacyclic group of order $2mn,$ where $m>2$ is odd. Then
\begin{enumerate}
\item[{\rm(a)}] $0\in\sigma(D^L(\Gamma_{M_{2mn}}))$ with multiplicity $1;$
\item[{\rm(b)}] $n(2m-1)\in\sigma(D^L(\Gamma_{M_{2mn}}))$  with multiplicity $m;$
\item[{\rm(c)}] $2mn\in\sigma(D^L(\Gamma_{M_{2mn}}))$  with multiplicity $m(n-1);$
\item[{\rm(d)}] $(3m-2)n\in\sigma(D^L(\Gamma_{M_{2mn}}))$  with multiplicity $(m-1)n-1.$
\end{enumerate}\et

\bt\label{dg4elt1} Let $\Gamma_{M_{2mn}}$ be the non commuting graph of the Metacyclic group of order $2mn,$ where $m>2$ is even. Then
\begin{enumerate}
\item[{\rm(a)}] $0\in\sigma(D^L(\Gamma_{M_{2mn}}))$ with multiplicity $1;$
\item[{\rm(b)}] $2n(m-1)\in\sigma(D^L(\Gamma_{M_{2mn}}))$  with multiplicity $\frac{m}{2};$
\item[{\rm(c)}] $2mn\in\sigma(D^L(\Gamma_{M_{2mn}}))$  with multiplicity $\frac{(2n-1)m}{2};$
\item[{\rm(d)}] $(3m-4)n\in\sigma(D^L(\Gamma_{M_{2mn}}))$  with multiplicity $(\frac{m}{2}-1)2n-1.$
\end{enumerate}\et

Hence, from Theorem~\ref{dg4olt1} and Theorem~\ref{dg4elt1}, it follows that $\Gamma_{M_{2mn}}$ is distance Laplacian integral for all $m$ and $n.$
\bt\label{dg4oqt1} Let $\Gamma_{M_{2mn}}$ be the non commuting graph of the Metacyclic group of order $2mn,$ where $m>2$ is odd. Then
\begin{enumerate}
\item[{\rm(a)}] $2mn-4\in\sigma(D^Q(\Gamma_{M_{2mn}}))$ with multiplicity $m(n-1);$
\item[{\rm(b)}] $(2m+1)n-4\in\sigma(D^Q(\Gamma_{M_{2mn}}))$ with multiplicity $m-1;$

\item[{\rm (c)}] $(3m-2)n-4 \in\sigma(D^Q(\Gamma_{M_{2mn}}))$ with multiplicity $(m-1)n-1;$

\item[{\rm (d)}] $nt_k(m-1)+(3mn+n-4) \in\sigma(D^Q(\Gamma_{M_{2mn}}))$ with multiplicity 1, where $t_k$ is a root of the equation $x^2(m-1)-(2m-5)t_k-m=0,$ for each $k=1,2.$

\end{enumerate}\et

Hence, from Theorem~\ref{dg4oqt1}, it follows that the distance signless Laplacian spectrum of $\Gamma_{M_{2mn}}$ is integral for those values of $m$ which gives integral values to $t_k,$ for $k=1,2.$

The case when $m$ is even is divided into two subcases. If $m=4,$  then $\Gamma_{M_{2mn}}=K_{2n,2n,2n}$ and we have the following result.
\bt\label{dg4eqtv1} Let $\Gamma_{M_{8n}}$ be the non commuting graph of the Metacyclic group of order $8n.$ Then
\begin{enumerate}
\item[{\rm(a)}] $8n-4\in\sigma(D^Q(\Gamma_{M_{8n}}))$ with multiplicity $3(2n-1);$
\item[{\rm(b)}] $10n-4\in\sigma(D^Q(\Gamma_{M_{8n}}))$ with multiplicity $2;$

\item[{\rm (c)}] $16n-4 \in\sigma(D^Q(\Gamma_{M_{8n}}))$ with multiplicity $1.$

\end{enumerate}\et

Hence, from Theorem~\ref{dg4eqtv1}, it follows that the distance signless Laplacian spectrum of $\Gamma_{M_{8n}}$ is integral for all values of $n.$ For any other even $m>4,$ the following result holds.
\bt\label{dg4eqtv2} Let $\Gamma_{M_{2mn}}$ be the non commuting graph of the Metacyclic group of order $2mn,$ where $m>4$ is even. Then
\begin{enumerate}
\item[{\rm(a)}] $3mn-4n-4\in\sigma(D^Q(\Gamma_{M_{2mn}}))$ with multiplicity $(m-2)n-1;$
\item[{\rm(b)}] $4mn-4n-4\in\sigma(D^Q(\Gamma_{M_{2mn}}))$ with multiplicity $\frac{(2n-1)m}{2};$

\item[{\rm (c)}] $2mn-4 \in\sigma(D^Q(\Gamma_{M_{2mn}}))$ with multiplicity $\frac{m}{2}-1;$

\item[{\rm (d)}] $nt_k(m-2)+(3mn+2n-4) \in\sigma(D^Q(\Gamma_{M_{2mn}}))$ with multiplicity 1, where $t_k$ is a root of the equation $x^2(m-2)-(m-5)2x-m=0,$ for each $k=1,2.$
 \end{enumerate}
\et

Hence, from Theorem~\ref{dg4eqtv2}, it follows that the distance signless Laplacian spectrum of $\Gamma_{M_{2mn}}$ is integral for those values of $m$ which gives integral values to $t_k,$ for $k=1,2.$

\section{Conclusion}
In this article, we have investigated the distance, distance (signless) Laplacian spectra of non-commuting graphs of the generalized quaternion group, the quasidihedral Group, $U_{6n},$ the metacyclic group, and obtain conditions under which these graphs will be distance, distance (signless) Laplacian integral.

\end{document}